\documentclass[12pt,twoside]{amsart}
\usepackage[latin1]{inputenc}
\usepackage{amsmath, amsthm, amscd, amsfonts, amssymb, graphicx}
\usepackage[bookmarksnumbered, plainpages]{hyperref}

\newtheorem{thm}{Theorem}[section]

\newtheorem{lem}[thm]{Lemma}
\newtheorem{prop}[thm]{Proposition}
\newtheorem{defn}[thm]{Definition}
\newtheorem{rem}[thm]{\bf{Remark}}

\numberwithin{equation}{section}

\begin{document}

\title{\bf The bimetric spectral Einstein-Hilbert action and the Kastler-Kalau-Walze type theorem for Lorentz warped products }
\author{Siyao Liu \hskip 0.4 true cm  Yong Wang$^{*}$}

\thanks{{\scriptsize
\hskip -0.4 true cm \textit{2010 Mathematics Subject Classification:}
53C40; 53C42.
\newline \textit{$^{*}$Corresponding author}}}

\maketitle

\begin{abstract}
 \indent In this paper, we define the bimetric spectral Einstein-Hilbert action which generalizes the spectral Einstein-Hilbert action.
 We compute the bimetric spectral Einstein-Hilbert action for the Lorentz warped product.
 Thus, we get the Kastler-Kalau-Walze type theorem for the Lorentz warped product.
\end{abstract}

\vskip 0.2 true cm


\pagestyle{myheadings}
\markboth{\rightline {\scriptsize Liu}}
         {\leftline{\scriptsize  The bimetric spectral Einstein-Hilbert action and the Kastler-Kalau-Walze type theorem for Lorentz warped products }}

\bigskip
\bigskip


\section{ Introduction }

Based on the noncommutative residue found in \cite{Gu,Wo}, Connes discovered the Kastler-Kalau-Walze type theorem as follow:
\begin{align}
Wres(D^{-2})=c_{0}\int_{M}sd{\rm Vol_{M} },\nonumber
\end{align}
where $M$ is an oriented spin Riemannian manifold, $D$ is the associated Dirac operator on the spinor bundle $S(TM),$ ${\rm Wres}$ denote the noncommutative residue, $s$ be the scalar curvature and $c_{0}$ be a constant.
The innovation of this theorem is to bring together two areas in which the noncommutative residue of the Dirac operator and the gravitational action (Einstein-Hilbert action).
This theorem gives a spectral explanation of the gravitational action.
The Kastler-Kalau-Walze type theorem was studied extensively by geometers \cite{Co1,Co2,K1,KW,Ac,U}.
Wang generalized some results to the case of manifolds with boundary in \cite{Wa1, Wa2} and proved the Kastler-Kalau-Walze type theorems for the Dirac operator and the signature operator on lower-dimensional manifolds with boundary.

When $(M,g)$ is a Lorentzian manifold, the Laplacian operator associated to $(M,g)$ is not a elliptic operator, so $\triangle^{-1}$ does not exist.
It is a natural question to get the Kastler-Kalau-Walze type theorem in the Lorentzian case.
In \cite{DW}, a Lorentz version of the Kastler-Kalau-Walze type theorem was given.
The concept of warped products was first introduced by Bishop and O'Neill (see \cite{BO}) to construct examples of Riemannian manifolds with negative curvature. 
In Riemannian geometry, warped product manifolds and their generic forms have been used to construct new examples with interesting curvature properties since then. 
{\bf The motivation of this paper} prove another Kastler-Kalau-Walze type theorem for the Lorentz warped product.

Our method is to introduce a spectral bimetric Einstein-Hilbert action $Wres(\triangle_{g_{1}}\triangle_{g_{2}}^{-2m})$ for metrics $g_{1}$ and $g_{2},$ $g_{2}$ is Riemannian metric and $g_{1}$ is any metric (including the Lorentz metric).
When $g_{1}=g_{2},$ we get the Kastler-Kalau-Walze type theorem.
When $g_{1}$ is Lorentz, we get a Lorentz Kastler-Kalau-Walze type theorem.
For this paper, we choose a classical space, namely the singly warped product, to investigate its generalization.
For the generalization of the singly warped product, at this point, the metric is not positive definite, and therefore the Laplace operator associated with it is not an elliptic differential operator.
The Laplace operator on the warped product is defined and the Kastler-Kalau-Walze type theorem for this non-elliptic differential operator is given computationally.
Our main theorem is as follows.

\begin{thm}\label{thm1.1}
Let $M$ be an $m$-dimensional compact Riemannian manifold and $N$ be an $n$-dimensional compact Riemannian manifold, for the definition of $g^{\varepsilon, f}$ and $g$ in Section 2, we get the following equality:
\begin{align}
&Wres [(\triangle^{_{\varepsilon}M\times_{f} N, g^{\varepsilon, f}})\circ(\triangle^{M\times N, g})^{-\overline{m}}]\nonumber\\
&=\int_{M\times N}\frac{2\pi^{\overline{m}}}{\Gamma(\overline{m})}\Big[\big(\frac{m-2}{12\varepsilon}+\frac{n}{12f^{2}}\big)S_{M}+\big(\frac{m}{12\varepsilon}+\frac{n-2}{12f^{2}}\big)S_{N}\Big]d{\rm Vol_{M\times N} }.\nonumber
\end{align}
where $2\overline{m}=m+n,$ $S_{M}$ and $S_{N}$ are scalar curvature of $(M, g^{M})$ and $(N, g^{N}).$
\end{thm}

\begin{rem}
In fact, Theorem 1.1 can be generalized for any manifolds without boundary. 
When $\varepsilon=1$ and $f=1,$ we get the Kastler-Kalau-Walze type theorem for $M\times N.$
$Wres [(\triangle^{_{\varepsilon}M\times_{f} N, g^{\varepsilon, f}})\circ(\triangle^{M\times N, g})^{-\overline{m}}]$ is the noncommutative integral of $\triangle^{_{\varepsilon}M\times_{f} N, g^{\varepsilon, f}}$ on $(M\times N, g)$ by \cite{WW}.
\end{rem}

A brief description of the organization of this paper is as follows.
In Section 2, this paper will first introduce the basic notions of the warped product.
In the next section, we recall some basic facts and formulas about the Laplace operator and also define the generalized Laplace operator.
The fourth section gives the research process of Theorem 1.1, which leads to the relevant conclusions for non-elliptic differential operators by classifying the computations. 
Give an example of the Kastler-Kalau-Walze type theorem for the Lorentz warped product in Section 5.

\vskip 1 true cm

\section{ The warped product }

\begin{defn}{\rm\cite{QW,W}}
The singly warped product $M\times_{f}N$ of two pseudo-Riemannian manifolds $(M, g^{M})$ and $(N, g^{N})$ with a smooth function $f: M \rightarrow (0, \infty)$ is a product manifold of form $M\times N$ with the metric tensor $g^{f}=g^{M}\oplus f^{2}g^{N}.$ 
\end{defn}
Here, $(M, g^{M})$ is called the base manifold and $(N, g^{N})$ is called as the fiber manifold and $f$ is called as the warping function.

A warped product is given as follows:
\begin{defn}
The warped product $_{\varepsilon}M\times_{f} N$ with $\varepsilon\neq0$ and a smooth function $f: M\rightarrow (0, \infty)$ for which $f>0$ is a product manifold of form $M\times N$ with the metric tensor $g^{\varepsilon, f}=\varepsilon g^{M}\oplus f^{2}g^{N}.$ 
\end{defn}

Based on the Koszul formula in {\rm\cite{K}}, we calculated that
\begin{prop}
Let $_{\varepsilon}M\times_{f} N$ be a warped product, $\nabla^{\varepsilon, f}, \nabla^{M}$ and $\nabla^{N}$ denote the Levi-Civita connection on $_{\varepsilon}M\times_{f} N, M$ and $N,$ respectively.
If $X, Y\in \Gamma(TM)$ and $U, V\in \Gamma(TN),$ then\\
(1)$\nabla^{\varepsilon,f}_{X}Y=\nabla^{M}_{X}Y;$\\
(2)$\nabla^{\varepsilon,f}_{X}U=\nabla^{\varepsilon,f}_{U}X=\frac{X(f)}{f}U;$\\
(3)$\nabla^{\varepsilon,f}_{U}V=-\frac{f}{\varepsilon}grad_{M}f g^{N}(U, V)+\nabla^{N}_{U}V.$
\end{prop}

\section{ The generalized Laplace operator }

Let $M$ be an $m$-dimensional compact Riemannian manifold with a Riemannian metric $g^{M}$ and $N$ be an $n$-dimensional compact Riemannian manifold with a Riemannian metric $g^{N}$.
The metric tensor of the product manifold $M\times N$ denoted as $g=g^{M}\oplus g^{N},$ the Levi-Civita connection of the product manifold $M\times N$ denoted as $\nabla.$
Set $\partial_{x_{i}}$ is a natural local frame on $TM,$ $(g^{ij}_{M})_{1\leq i, j\leq m}$ is the inverse matrix associated to the metric matrix $(g^{M}_{ij})_{1\leq i, j\leq m}$ on $M;$ $\partial_{x_{\alpha}}$ is a the natural local frame on $TN,$ $(g^{\alpha\beta}_{N})_{m+1\leq \alpha, \beta\leq m+n}$ is the inverse matrix associated to the metric matrix $(g^{N}_{\alpha\beta})_{m+1\leq \alpha, \beta\leq m+n}$ on $N.$ 

According to the definition in {\rm\cite{DSZ}}, we can write the Laplace operator on $M\times N$ as follows:
\begin{align}
\triangle^{M\times N, g}=-\sum_{j, l=1}^{m+n}g^{jl}(\partial_{j}\partial_{l}-\sum_{\alpha=1}^{m+n}\Gamma_{ij}^{\alpha}\partial_{\alpha}),\nonumber
\end{align}
where $\Gamma_{ij}^{\alpha}$ is the Christoffel coefficient of Levi-Civita connection $\nabla.$

It is immediate that
\begin{align}
[g_{\varepsilon,f}]= \left[\begin{array}{lcr}
  \varepsilon[g_{ij}^{M}]  & 0  \\
   0  &  f^{2}[g_{\alpha\beta}^{N}]
\end{array}\right];~~~
[(g^{-1})_{\varepsilon,f}]=[g^{\varepsilon,f}]= \left[\begin{array}{lcr}
  \frac{1}{\varepsilon}[g^{ij}_{M}]  & 0  \\
   0  &  \frac{1}{f^{2}}[g^{\alpha\beta}_{N}]
\end{array}\right],\nonumber
\end{align}
where $1\leq i,j\leq m, m+1\leq \alpha,\beta\leq m+n.$

Similar to the definition of $\triangle^{M\times N, g},$ we have the generalized Laplace operator on $_{\varepsilon}M\times_{f} N$ that is
\begin{align}
\triangle^{_{\varepsilon}M\times_{f} N, g^{\varepsilon, f}}=-\sum_{j, l=1}^{m+n}g^{jl}_{\varepsilon, f}(\partial_{j}\partial_{l}-\sum_{\alpha=1}^{m+n}\Gamma_{ij}^{\alpha;\varepsilon, f}\partial_{\alpha}),\nonumber
\end{align}
where $\Gamma_{ij}^{\alpha;\varepsilon, f}$ is the Christoffel coefficient of Levi-Civita connection $\nabla^{\varepsilon, f},$ see {\rm\cite{O}} for more details.

\section{ The noncommutative integral }
Extending previous definitions by Connes {\rm\cite{C}}, a noncommutative integral was introduced in {\rm\cite{FGL}} based on the noncommutative residue {\rm\cite{W1}}, combine (1.4) in {\rm\cite{CC}} and {\rm\cite{K1}}, using the definition of the residue, gives the noncommutative residue of the pseudo-differential operator $P$ and the operator $(\triangle^{M\times N, g})^{-\overline{m}}$:
\begin{align}
\int\hspace{-1.05em}-\ P\ ds^{2\overline{m}}:=Wres[P(\triangle^{M\times N, g})^{-\overline{m}}]:=\int_{M\times N}\int_{S^{*}(M\times N)}{\rm
tr}\big[\sigma_{-2\overline{m}}\big(P(\triangle^{M\times N, g})^{-\overline{m}}\big)\big](x, \xi),\nonumber
\end{align}
where ${\rm tr}$ as shorthand of ${\rm trace}.$

The aim of this section is to calculate
\begin{align}
&\int\hspace{-1.05em}-\ \triangle^{_{\varepsilon}M\times_{f} N, g^{\varepsilon, f}}\ ds^{2\overline{m}}=Wres [(\triangle^{_{\varepsilon}M\times_{f} N, g^{\varepsilon, f}})\circ(\triangle^{M\times N, g})^{-\overline{m}}]\\
&=\int_{M\times N}\int_{S^{*}(M\times N)}{\rm
tr}\big[\sigma_{-2\overline{m}}\big((\triangle^{_{\varepsilon}M\times_{f} N, g^{\varepsilon, f}})\circ(\triangle^{M\times N, g})^{-\overline{m}}\big)\big](x, \xi).\nonumber
\end{align}

Suppose $\varepsilon=1$ and $f=1,$ then Equation (4.1) becomes the Kastler-Kalau-Walze type theorem, implying that the noncommutative integral under certain conditions is the Kastler-Kalau-Walze type theorem.
When $\varepsilon=-1,$ we get the Lorentz Kastler-Kalau-Walze type theorem for $M\times N.$

From \cite{W2}, we similarly can obtain
\begin{align}
&\sigma_{-2\overline{m}}\big((\triangle^{_{\varepsilon}M\times_{f} N, g^{\varepsilon, f}})\circ(\triangle^{M\times N, g})^{-\overline{m}}\big)\nonumber\\
&=\Bigg\{\sum_{|\alpha|=0}^{+\infty}\frac{1}{\alpha!}\partial_{\xi}^{\alpha}[\sigma(\triangle^{_{\varepsilon}M\times_{f} N, g^{\varepsilon, f}})]D_{x}^{\alpha}\big[\sigma\big((\triangle^{M\times N, g})^{-\overline{m}}\big)\big]\Bigg\}_{-2\overline{m}}\nonumber\\
&=\sigma_{0}(\triangle^{_{\varepsilon}M\times_{f} N, g^{\varepsilon, f}})\sigma_{-2\overline{m}}\big((\triangle^{M\times N, g})^{-\overline{m}}\big)+\sigma_{1}(\triangle^{_{\varepsilon}M\times_{f} N, g^{\varepsilon, f}})\sigma_{-2\overline{m}-1}\big((\triangle^{M\times N, g})^{-\overline{m}}\big)\nonumber\\
&+\sigma_{2}(\triangle^{_{\varepsilon}M\times_{f} N, g^{\varepsilon, f}})\sigma_{-2\overline{m}-2}\big((\triangle^{M\times N, g})^{-\overline{m}}\big)-i\sum_{\lambda=1}^{2\overline{m}}\partial_{\xi_{\lambda}}[\sigma_{2}(\triangle^{_{\varepsilon}M\times_{f} N, g^{\varepsilon, f}})]\nonumber\\
&\cdot\partial_{x_{\lambda}}\big[\sigma_{-2\overline{m}-1}\big((\triangle^{M\times N, g})^{-\overline{m}}\big)\big]-i\sum_{\lambda=1}^{2\overline{m}}\partial_{\xi_{\lambda}}[\sigma_{1}(\triangle^{_{\varepsilon}M\times_{f} N, g^{\varepsilon, f}})]\partial_{x_{\lambda}}\big[\sigma_{-2\overline{m}}\big((\triangle^{M\times N, g})^{-\overline{m}}\big)\big]\nonumber\\
&-\frac{1}{2}\sum_{\lambda,\nu=1}^{2\overline{m}}\partial_{\xi_{\lambda}}\partial_{\xi_{\nu}}[\sigma_{2}(\triangle^{_{\varepsilon}M\times_{f} N, g^{\varepsilon, f}})]\partial_{x_{\lambda}}\partial_{x_{\nu}}\big[\sigma_{-2\overline{m}}\big((\triangle^{M\times N, g})^{-\overline{m}}\big)\big].\nonumber
\end{align}

A simple calculation shows that
\begin{lem}
The symbols of $\triangle^{M\times N, g}$ are given
\begin{align}
\sigma_{2}(\triangle^{M\times N, g})&=|\xi|^{2}_{g};\nonumber\\
\sigma_{1}(\triangle^{M\times N, g})&=\sqrt{-1}\sum_{j, l, \alpha=1}^{m+n}g^{jl}\Gamma_{ij}^{\alpha}\xi_{\alpha};\nonumber\\
\sigma_{0}(\triangle^{M\times N, g})&=0.\nonumber
\end{align}
\end{lem}

Write
 \begin{align}
D_x^{\alpha}&=(-i)^{|\alpha|}\partial_x^{\alpha};
~\sigma(\triangle^{M\times N, g})=p_{2}+p_{1}+p_{0};
~\sigma\big((\triangle^{M\times N, g})^{-1}\big)=\sum^{\infty}_{j=2}q_{-j}.\nonumber
\end{align}
By the composition formula of pseudodifferential operators, we have
\begin{align}
1&=\sigma[\triangle^{M\times N, g}\circ (\triangle^{M\times N, g})^{-1}]=
\sum_{\alpha}\frac{1}{\alpha!}\partial^{\alpha}_{\xi}
[\sigma(\triangle^{M\times N, g})]{{D}}^{\alpha}_{x}
\big[\sigma\big((\triangle^{M\times N, g})^{-1}\big)\big]\nonumber\\
&=(p_2+p_1+p_0)(q_{-2}+q_{-3}+\cdots)\nonumber\\
&+\sum_{j}(\partial_{\xi_{j}}p_2+\partial_{\xi_j}p_1+\partial_{\xi_{j}}p_0)
(D_{x_j}q_{-2}+D_{x_j}q_{-3}+\cdots)\nonumber\\
&+\frac{1}{2}\sum_{i,j}(\partial_{\xi_{i}}\partial_{\xi_{j}}p_2+\partial_{\xi_{i}}\partial_{\xi_j}p_1+\partial_{\xi_{i}}\partial_{\xi_{j}}p_0)
(D_{x_i}D_{x_j}q_{-2}+D_{x_i}D_{x_j}q_{-3}+\cdots) \nonumber\\
&=p_2q_{-2}+(p_2q_{-3}+p_1q_{-2}+\sum_j\partial_{\xi_j}p_2D_{x_j}q_{-2})+\cdots,\nonumber
\end{align}
so
\begin{align}
q_{-2}&=p_2^{-1};\nonumber\\
q_{-3}&=-p_2^{-1}[p_1q_{-2}+\sum_j\partial_{\xi_j}p_2D_{x_j}(q_{-2})];\nonumber\\
q_{-4}&=-p_2^{-1}[p_1q_{-3}+p_0q_{-2}+\sum_j\partial_{\xi_j}p_1D_{x_j}(q_{-2})+\sum_j\partial_{\xi_j}p_2D_{x_j}(q_{-3})\\
&+\frac{1}{2}\sum_{i,j}\partial_{\xi_i}\partial_{\xi_j}p_2D_{x_i}D_{x_j}(q_{-2})].\nonumber
\end{align}
Then 
\begin{lem} The following identities hold:
\begin{align}
\sigma_{-2}\big((\triangle^{M\times N, g})^{-1}\big)&=|\xi|_{g}^{-2};\nonumber\\
\sigma_{-3}\big((\triangle^{M\times N, g})^{-1}\big)&=-\sqrt{-1}|\xi|_{g}^{-4}\sum^{m+n}_{k=1}\Gamma^{k}\xi_{k}-2\sqrt{-1}|\xi|_{g}^{-6}\sum^{m+n}_{j,\alpha,\beta=1}\xi^{j}\xi_{\alpha}\xi_{\beta}\partial_{x_{j}}(g^{\alpha\beta}).\nonumber
\end{align}
\end{lem}

Since $\sigma_{-2}\big((\triangle^{M\times N, g})^{-1}\big)=|\xi|_{g}^{-2},$ it follows that 
\begin{align}
\sigma_{-2\overline{m}}\big((\triangle^{M\times N, g})^{-\overline{m}}\big)&=|\xi|_{g}^{-2\overline{m}}.\nonumber
\end{align}

From \cite{W2}, we similarly can obtain
\begin{align}
&\sigma_{-2\overline{m}-1}\big((\triangle^{M\times N, g})^{-\overline{m}}\big)=\overline{m}\sigma_{2}(\triangle^{M\times N, g})^{(-\overline{m}+1)}\sigma_{-3}\big((\triangle^{M\times N, g})^{-1}\big)\nonumber\\
&-\sqrt{-1}\sum_{\mu=0}^{2\overline{m}}\sum_{k=0}^{\overline{m}-2}\partial_{\xi_{\mu}}[\sigma_{2}(\triangle^{M\times N, g})^{(-\overline{m}+k+1)}]\partial_{x_{\mu}}[\sigma_{2}(\triangle^{M\times N, g})^{(-1)}]\big(\sigma_{2}(\triangle^{M\times N, g})\big)^{(-k)}\nonumber\\
&=\overline{m}|\xi|_{g}^{(-2\overline{m}+2)}[-\sqrt{-1}|\xi|_{g}^{-4}\sum_{t=1}^{2\overline{m}}\Gamma^{t}\xi_{t}-2\sqrt{-1}|\xi|_{g}^{-6}\sum^{2\overline{m}}_{j,\alpha,\beta=1}\xi^{j}\xi_{\alpha}\xi_{\beta}\partial_{x_{j}}(g^{\alpha\beta})]\nonumber\\
&+2\sqrt{-1}\sum_{\mu=1}^{2\overline{m}}\sum_{k=0}^{\overline{m}-2}(-\overline{m}+k+1)|\xi|_{g}^{(-2\overline{m}-4)}\xi^{\mu}\xi_{\alpha}\xi_{\beta}\partial_{x_{\mu}}(g^{\alpha\beta}).\nonumber
\end{align}

Combined with \cite{KW} and \cite{WW}, we get the recursion relations
\begin{align}
&\sigma_{-2\overline{m}-2}\big((\triangle^{M\times N, g})^{-\overline{m}}\big)(x_{0})\nonumber\\
&=\sum_{|\alpha|=0}^{2}\sum_{k=2}^{4-|\alpha|}(-\sqrt{-1})^{|\alpha|}\frac{1}{\alpha!}\partial_{\xi}^{\alpha}\big[\sigma_{|\alpha|+k-2\overline{m}-2}\big((\triangle^{M\times N, g})^{-\overline{m}+1}\big)(x_{0})\big]\partial_{x}^{\alpha}\big[\sigma_{-k}\big((\triangle^{M\times N, g})^{-1}\big)(x_{0})\big]\nonumber\\
&=\frac{\overline{m}}{4}\big[(2\overline{m}+2)|\xi|_{g}^{(-2\overline{m}+2)}\sigma_{-4}\big((\triangle^{M\times N, g})^{-1}\big)(x_{0})+(2\overline{m}-2)|\xi|_{g}^{(-2\overline{m}-2)}\sigma_{0}(\triangle^{M\times N, g})(x_{0})\big]\nonumber\\
&=\frac{\overline{m}(\overline{m}+1)}{3}|\xi|_{g}^{(-2\overline{m}-4)}\sum_{j,l,\gamma=1}^{m+n}\xi_{j}\xi_{l}R_{\gamma j\gamma l}(x_{0}).\nonumber
\end{align}

Furthermore,
\begin{lem}
The symbols of $\triangle^{_{\varepsilon}M\times_{f} N, g^{\varepsilon, f}}$ are given
\begin{align}
\sigma_{2}(\triangle^{_{\varepsilon}M\times_{f} N, g^{\varepsilon, f}})&=\frac{1}{\varepsilon}|\xi^{M}|^{2}_{g^{M}}+\frac{1}{f^{2}}|\xi^{N}|^{2}_{g^{N}};\nonumber\\
\sigma_{1}(\triangle^{_{\varepsilon}M\times_{f} N, g^{\varepsilon, f}})&=\frac{\sqrt{-1}}{\varepsilon}\sum_{j, l, k=1}^{m}g^{jl}_{M}\Gamma_{jl}^{k;M}\xi_{k}-\frac{\sqrt{-1}n}{\varepsilon f}\sum_{j, k=1}^{m}\frac{\partial f}{\partial x_{j}}g^{jk}_{M}\xi_{k}\nonumber\\
&+\frac{\sqrt{-1}}{f^{2}}\sum_{\alpha, \beta, \gamma=m+1}^{m+n}g^{\alpha\beta}_{N}\Gamma_{\alpha\beta}^{\gamma;N}\xi_{\gamma};\nonumber\\
\sigma_{0}(\triangle^{_{\varepsilon}M\times_{f} N, g^{\varepsilon, f}})&=0.\nonumber
\end{align}
\end{lem}

In the next step, recall the technical tool used in the calculation which is the integration of polynomial functions on the unit sphere.
According to \cite{WW}, we have 
\begin{align}
I^{\gamma_{1}\cdot\cdot\cdot\gamma_{2\overline{n}+2}}_{S_{n}}=\int_{|x|=1}d^{n}xx^{\gamma_{1}}\cdot\cdot\cdot x^{\gamma_{2\overline{n}+2}},\nonumber
\end{align}
namely, the monomial integral on the unit sphere.
By Proposition A.2 in \cite{BS}, the polynomial integrals on the $n$-dimesional sphere are given
\begin{align}
I^{\gamma_{1}\cdot\cdot\cdot\gamma_{2\overline{n}+2}}_{S_{n}}=\frac{1}{2\overline{n}+n}[\delta^{\gamma_{1}\gamma_{2}}I^{\gamma_{3}\cdot\cdot\cdot\gamma_{2\overline{n}+2}}_{S_{n}}+\cdot\cdot\cdot+\delta^{\gamma_{1}\gamma_{2\overline{n}+1}}I^{\gamma_{2}\cdot\cdot\cdot\gamma_{2\overline{n}+1}}_{S_{n}}],\nonumber
\end{align}
where $S_{n}=S^{n-1}$ in $\mathbb{R}^{n}.$
For $\overline{n}=0,$ we have $I^{0}=area(S_{n})=\frac{2\pi^{\frac{n}{2}}}{\Gamma(\frac{n}{2})},$ we immediately get
\begin{align}
&I^{\gamma_{1}\gamma_{2}}_{S_{n}}=\frac{1}{n}area(S_{n})\delta^{\gamma_{1}\gamma_{2}};\nonumber
\end{align}
\begin{align}
&I^{\gamma_{1}\gamma_{2}\gamma_{3}\gamma_{4}}_{S_{n}}=\frac{1}{n(n+2)}area(S_{n})[\delta^{\gamma_{1}\gamma_{2}}\delta^{\gamma_{3}\gamma_{4}}+\delta^{\gamma_{1}\gamma_{3}}\delta^{\gamma_{2}\gamma_{4}}+\delta^{\gamma_{1}\gamma_{4}}\delta^{\gamma_{2}\gamma_{3}}].\nonumber
\end{align}

Beside, let us start to compute each term of $\int_{S^{*}(M\times N)}{\rm
tr}\big[\sigma_{-2\overline{m}}\big((\triangle^{_{\varepsilon}M\times_{f} N, g^{\varepsilon, f}})\circ(\triangle^{M\times N, g})^{-\overline{m}}\big)\big]$
$(x, \xi).$

\noindent {\bf $\uppercase\expandafter{\romannumeral1}$)}~For $\sigma_{0}(\triangle^{_{\varepsilon}M\times_{f} N, g^{\varepsilon, f}})\sigma_{-2\overline{m}}\big((\triangle^{M\times N, g})^{-\overline{m}}\big):$\\

We have $\sigma_{0}(\triangle^{_{\varepsilon}M\times_{f} N, g^{\varepsilon, f}})=0,$ so
\begin{align}
\int_{|\xi|=1}{\rm tr}\big[\sigma_{0}(\triangle^{_{\varepsilon}M\times_{f} N, g^{\varepsilon, f}})\sigma_{-2\overline{m}}\big((\triangle^{M\times N, g})^{-\overline{m}}\big)(x_{0})\big]\sigma(\xi)=0.\nonumber
\end{align}

\noindent {\bf $\uppercase\expandafter{\romannumeral2}$)}~For $\sigma_{1}(\triangle^{_{\varepsilon}M\times_{f} N, g^{\varepsilon, f}})\sigma_{-2\overline{m}-1}\big((\triangle^{M\times N, g})^{-\overline{m}}\big):$\\

Let $x_{0}=(\overline{x}_{0},\widetilde{x}_{0}),$ we take the normal coordinates about $\overline{x}_{0}$ in $(M, g^{M})$ and $\widetilde{x}_{0}$ in $(N, g^{N}).$

Using the facts
\begin{align}
\Gamma^{\mu}(x_{0})=\partial_{x_{\mu}}(g^{\alpha\beta})(x_{0})=0,~\partial_{x_{k}}(\Gamma^{\mu})(x_{0})=\frac{2}{3}R_{k\alpha\mu\alpha}(x_{0}),\nonumber
\end{align}
we see that
\begin{align}
\sigma_{-2\overline{m}-1}\big((\triangle^{M\times N, g})^{-\overline{m}}\big)(x_{0})=0.\nonumber
\end{align}

Thus,
\begin{align}
\int_{|\xi|=1}{\rm tr}\big[\sigma_{1}(\triangle^{_{\varepsilon}M\times_{f} N, g^{\varepsilon, f}})\sigma_{-2\overline{m}-1}\big((\triangle^{M\times N, g})^{-\overline{m}}\big)(x_{0})\big]\sigma(\xi)=0.\nonumber
\end{align}

\noindent {\bf $\uppercase\expandafter{\romannumeral3}$)}~For $\sigma_{2}(\triangle^{_{\varepsilon}M\times_{f} N, g^{\varepsilon, f}})\sigma_{-2\overline{m}-2}\big((\triangle^{M\times N, g})^{-\overline{m}}\big):$\\

In normal coordinates, (4.2) now becomes
\begin{align}
\sigma_{-4}\big((\triangle^{M\times N, g})^{-1}\big)(x_{0})=q_{-4}(x_{0})=\frac{2}{3}|\xi|_{g}^{-6}\sum_{j,l,\gamma=1}^{m+n}\xi_{j}\xi_{l}R_{\gamma j\gamma l}(x_{0}).\nonumber
\end{align}

Meanwhile, we can obtain
\begin{align}
\sigma_{2}(\triangle^{_{\varepsilon}M\times_{f} N, g^{\varepsilon, f}})\sigma_{-2\overline{m}-2}\big((\triangle^{M\times N, g})^{-\overline{m}}\big)(x_{0})&=\frac{\overline{m}(\overline{m}+1)}{3\varepsilon}|\xi|_{g}^{(-2\overline{m}-4)}\sum_{j,l,\gamma=1}^{m+n}\sum_{a=1}^{m}\xi_{j}\xi_{l}\xi^{2}_{a}R_{\gamma j\gamma l}(x_{0})\nonumber\\
&+\frac{\overline{m}(\overline{m}+1)}{3f^{2}}|\xi|_{g}^{(-2\overline{m}-4)}\sum_{j,l,\gamma=1}^{m+n}\sum_{c=m+1}^{m+n}\xi_{j}\xi_{l}\xi^{2}_{c}R_{\gamma j\gamma l}(x_{0}),\nonumber
\end{align}
by taking into account
\begin{align}
\sum_{j,l,\gamma=1}^{m+n}\sum_{a=1}^{m}\int_{|\xi|=1}\xi_{j}\xi_{l}\xi^{2}_{a}\sigma(\xi)&=\frac{1}{2\overline{m}(2\overline{m}+2)}area(S_{2\overline{m}})\sum_{j,l,\gamma=1}^{m+n}\sum_{a=1}^{m}[\delta^{aa}\delta^{jl}+\delta^{aj}\delta^{al}+\delta^{al}\delta^{aj}]\nonumber\\
&=\frac{1}{2\overline{m}(2\overline{m}+2)}area(S_{2\overline{m}})\sum_{j,l,\gamma=1}^{m+n}\sum_{a=1}^{m}[m\delta^{jl}+\delta^{aj}\delta^{al}+\delta^{al}\delta^{aj}]\nonumber
\end{align}
and
\begin{align}
\sum_{j,l,\gamma=1}^{m+n}\sum_{c=m+1}^{m+n}\int_{|\xi|=1}\xi_{j}\xi_{l}\xi^{2}_{c}\sigma(\xi)&=\frac{1}{2\overline{m}(2\overline{m}+2)}area(S_{2\overline{m}})\sum_{j,l,\gamma=1}^{m+n}\sum_{c=m+1}^{m+n}[\delta^{cc}\delta^{jl}+\delta^{cj}\delta^{cl}+\delta^{cl}\delta^{cj}]\nonumber\\
&=\frac{1}{2\overline{m}(2\overline{m}+2)}area(S_{2\overline{m}})\sum_{j,l,\gamma=1}^{m+n}\sum_{c=m+1}^{m+n}[n\delta^{jl}+\delta^{cj}\delta^{cl}+\delta^{cl}\delta^{cj}]\nonumber,
\end{align}
we have
\begin{align}
&\int_{|\xi|=1}{\rm tr}\big[\sigma_{2}(\triangle^{_{\varepsilon}M\times_{f} N, g^{\varepsilon, f}})\sigma_{-2\overline{m}-2}\big((\triangle^{M\times N, g})^{-\overline{m}}\big)(x_{0})\big]\sigma(\xi)\nonumber\\
&=\frac{m}{12\varepsilon}area(S_{2\overline{m}})\sum_{j,\gamma=1}^{m+n}R_{\gamma j\gamma j}(x_{0})+\frac{1}{6\varepsilon}area(S_{2\overline{m}})\sum_{\gamma=1}^{m+n}\sum_{a=1}^{m}R_{\gamma a\gamma a}(x_{0})\nonumber\\
&+\frac{n}{12f^{2}}area(S_{2\overline{m}})\sum_{j,\gamma=1}^{m+n}R_{\gamma j\gamma j}(x_{0})+\frac{1}{6f^{2}}area(S_{2\overline{m}})\sum_{\gamma=1}^{m+n}\sum_{c=m+1}^{m+n}R_{\gamma c\gamma c}(x_{0}).\nonumber
\end{align}

\noindent {\bf $\uppercase\expandafter{\romannumeral4}$)}~For $-i\sum_{\lambda=1}^{2\overline{m}}\partial_{\xi_{\lambda}}[\sigma_{2}(\triangle^{_{\varepsilon}M\times_{f} N, g^{\varepsilon, f}})]\partial_{x_{\lambda}}\big[\sigma_{-2\overline{m}-1}\big((\triangle^{M\times N, g})^{-\overline{m}}\big)\big]:$\\

By using (2.5) in \cite{W1} and making tedious calculations, we have
\begin{align}
\partial_{x_{j}}\partial_{x_{l}}(g^{\alpha\beta})(x_{0})=\partial_{x_{j}}\partial_{x_{l}}\big(\frac{1}{3}R_{\alpha c\beta d}(x_{0})x_{c}x_{d}\big)=\frac{1}{3}R_{\alpha l\beta j}(x_{0})+\frac{1}{3}R_{\alpha j\beta l}(x_{0}).\nonumber
\end{align}

Computations show that
\begin{align}
\sum_{\lambda=1}^{2\overline{m}}\partial_{\xi_{\lambda}}[\sigma_{2}(\triangle^{_{\varepsilon}M\times_{f} N, g^{\varepsilon, f}})](x_{0})&=\sum_{\lambda=1}^{2\overline{m}}\partial_{\xi_{\lambda}}\big(\frac{1}{\varepsilon}\sum_{a,b=1}^{m}g^{ab}_{M}\xi_{a}\xi_{b}+\frac{1}{f^{2}}\sum_{c,d=m+1}^{m+n}g^{cd}_{N}\xi_{c}\xi_{d}\big)(x_{0})\nonumber\\
&=\frac{2}{\varepsilon}\sum_{\lambda=1}^{2\overline{m}}\sum_{a=1}^{m}\partial_{\xi_{\lambda}}(\xi_{a})\xi_{a}+\frac{2}{f^{2}}\sum_{\lambda=1}^{2\overline{m}}\sum_{c=m+1}^{m+n}\partial_{\xi_{\lambda}}(\xi_{c})\xi_{c}\nonumber
\end{align}
and
\begin{align}
&\sum_{\lambda=1}^{2\overline{m}}\partial_{x_{\lambda}}\big[\sigma_{-2\overline{m}-1}\big((\triangle^{M\times N, g})^{-\overline{m}}\big)\big](x_{0})\nonumber\\
&=-\overline{m}\sqrt{-1}|\xi|_{g}^{(-2\overline{m}-2)}\sum_{\lambda,t=1}^{2\overline{m}}\partial_{x_{\lambda}}(\Gamma^{t})(x_{0})\xi_{t}\nonumber\\
&-2\overline{m}\sqrt{-1}|\xi|_{g}^{(-2\overline{m}-4)}\sum^{2\overline{m}}_{\lambda,j,\alpha,\beta=1}\xi^{j}\xi_{\alpha}\xi_{\beta}\partial_{x_{\lambda}}\partial_{x_{j}}(g^{\alpha\beta})(x_{0})\nonumber\\
&+2\sqrt{-1}|\xi|_{g}^{(-2\overline{m}-4)}\sum_{\lambda,\mu,\alpha,\beta=1}^{2\overline{m}}\sum_{k=0}^{\overline{m}-2}(-\overline{m}+k+1)\xi^{\mu}\xi_{\alpha}\xi_{\beta}\partial_{x_{\lambda}}\partial_{x_{\mu}}(g^{\alpha\beta})(x_{0})\nonumber\\
&=-\frac{2\overline{m}\sqrt{-1}}{3}|\xi|_{g}^{(-2\overline{m}-2)}\sum_{\lambda,t,\alpha=1}^{2\overline{m}}R_{\alpha \lambda\alpha t}(x_{0})\xi_{t}.\nonumber
\end{align}

From this
\begin{align}
&-i\sum_{\lambda=1}^{2\overline{m}}\partial_{\xi_{\lambda}}[\sigma_{2}(\triangle^{_{\varepsilon}M\times_{f} N, g^{\varepsilon, f}})]\partial_{x_{\lambda}}\big[\sigma_{-2\overline{m}-1}\big((\triangle^{M\times N, g})^{-\overline{m}}\big)\big](x_{0})\nonumber\\
&=-\frac{4\overline{m}}{3\varepsilon}|\xi|_{g}^{(-2\overline{m}-2)}\sum_{\lambda,t,\alpha=1}^{2\overline{m}}\sum_{a=1}^{m}R_{\alpha \lambda\alpha t}(x_{0})\partial_{\xi_{\lambda}}(\xi_{a})\xi_{a}\xi_{t}\nonumber\\
&-\frac{4\overline{m}}{3f^{2}}|\xi|_{g}^{(-2\overline{m}-2)}\sum_{\lambda,t,\alpha=1}^{2\overline{m}}\sum_{c=m+1}^{m+n}R_{\alpha \lambda\alpha t}(x_{0})\partial_{\xi_{\lambda}}(\xi_{c})\xi_{c}\xi_{t},\nonumber
\end{align}
so
\begin{align}
&\int_{|\xi|=1}{\rm tr}\big[-i\sum_{\lambda=1}^{2\overline{m}}\partial_{\xi_{\lambda}}[\sigma_{2}(\triangle^{_{\varepsilon}M\times_{f} N, g^{\varepsilon, f}})]\partial_{x_{\lambda}}\big[\sigma_{-2\overline{m}-1}\big((\triangle^{M\times N, g})^{-\overline{m}}\big)\big](x_{0})\big]\sigma(\xi)\nonumber\\
&=-\frac{2}{3\varepsilon}area(S_{2\overline{m}})\sum_{\alpha=1}^{2\overline{m}}\sum_{a=1}^{m}R_{\alpha a\alpha a}(x_{0})-\frac{2}{3f^{2}}area(S_{2\overline{m}})\sum_{\alpha=1}^{2\overline{m}}\sum_{c=m+1}^{m+n}R_{\alpha c\alpha c}(x_{0}).\nonumber
\end{align}

\noindent {\bf $\uppercase\expandafter{\romannumeral5}$)}~For $-i\sum_{\lambda=1}^{2\overline{m}}\partial_{\xi_{\lambda}}[\sigma_{1}(\triangle^{_{\varepsilon}M\times_{f} N, g^{\varepsilon, f}})]\partial_{x_{\lambda}}\big[\sigma_{-2\overline{m}}\big((\triangle^{M\times N, g})^{-\overline{m}}\big)\big]:$\\

We deduce that 
\begin{align}
&\sum_{\lambda=1}^{2\overline{m}}\partial_{x_{\lambda}}\big[\sigma_{-2\overline{m}}\big((\triangle^{M\times N, g})^{-\overline{m}}\big)\big](x_{0})=\sum_{\lambda,\alpha,\beta=1}^{2\overline{m}}\xi_{\alpha}\xi_{\beta}\partial_{x_{\lambda}}(g^{\alpha\beta})(x_{0})=0.\nonumber
\end{align}

This gives 
\begin{align}
&\int_{|\xi|=1}{\rm tr}\big[-i\sum_{\lambda=1}^{2\overline{m}}\partial_{\xi_{\lambda}}[\sigma_{1}(\triangle^{_{\varepsilon}M\times_{f} N, g^{\varepsilon, f}})]\partial_{x_{\lambda}}\big[\sigma_{-2\overline{m}}\big((\triangle^{M\times N, g})^{-\overline{m}}\big)\big](x_{0})\big]\sigma(\xi)=0.\nonumber
\end{align}

\noindent {\bf $\uppercase\expandafter{\romannumeral6}$)}~For $-\frac{1}{2}\sum_{\lambda,\nu=1}^{2\overline{m}}\partial_{\xi_{\lambda}}\partial_{\xi_{\nu}}[\sigma_{2}(\triangle^{_{\varepsilon}M\times_{f} N, g^{\varepsilon, f}})]\partial_{x_{\lambda}}\partial_{x_{\nu}}\big[\sigma_{-2\overline{m}}\big((\triangle^{M\times N, g})^{-\overline{m}}\big)\big]:$\\

Let us first prove that
\begin{align}
&\sum_{\lambda,\nu=1}^{2\overline{m}}\partial_{\xi_{\lambda}}\partial_{\xi_{\nu}}[\sigma_{2}(\triangle^{_{\varepsilon}M\times_{f} N, g^{\varepsilon, f}})](x_{0})\nonumber\\
&=\frac{2}{\varepsilon}\sum_{\lambda,\nu=1}^{2\overline{m}}\sum_{a=1}^{m}\partial_{\xi_{\lambda}}(\xi_{a})\partial_{\xi_{\nu}}(\xi_{a})+\frac{2}{f^{2}}\sum_{\lambda,\nu=1}^{2\overline{m}}\sum_{c=m+1}^{m+n}\partial_{\xi_{\lambda}}(\xi_{c})\partial_{\xi_{\nu}}(\xi_{c}).\nonumber
\end{align}
Likewise,
\begin{align}
\sum_{\lambda,\nu=1}^{2\overline{m}}\partial_{x_{\lambda}}\partial_{x_{\nu}}\big[\sigma_{-2\overline{m}}\big((\triangle^{M\times N, g})^{-\overline{m}}\big)\big](x_{0})&=-\overline{m}|\xi|_{g}^{(-2\overline{m}-2)}\sum_{\lambda,\nu,\alpha,\beta=1}^{2\overline{m}}\xi_{\alpha}\xi_{\beta}\partial_{x_{\nu}}\partial_{x_{\lambda}}(g^{\alpha\beta})(x_{0})\nonumber\\
&=-\frac{2\overline{m}}{3}|\xi|_{g}^{(-2\overline{m}-2)}\sum_{\lambda,\nu,\alpha,\beta=1}^{2\overline{m}}R_{\alpha\lambda\beta\nu}(x_{0})\xi_{\alpha}\xi_{\beta}.\nonumber
\end{align}

The result is 
\begin{align}
&\int_{|\xi|=1}{\rm tr}\big[-\frac{1}{2}\sum_{\lambda,\nu=1}^{2\overline{m}}\partial_{\xi_{\lambda}}\partial_{\xi_{\nu}}[\sigma_{2}(\triangle^{_{\varepsilon}M\times_{f} N, g^{\varepsilon, f}})]\partial_{x_{\lambda}}\partial_{x_{\nu}}\big[\sigma_{-2\overline{m}}\big((\triangle^{M\times N, g})^{-\overline{m}}\big)\big](x_{0})\big]\sigma(\xi)\nonumber\\
&=\frac{1}{3\varepsilon}area(S_{2\overline{m}})\sum_{\alpha=1}^{2\overline{m}}\sum_{a=1}^{m}R_{\alpha a\alpha a}(x_{0})+\frac{1}{3f^{2}}area(S_{2\overline{m}})\sum_{\alpha=1}^{2\overline{m}}\sum_{c=m+1}^{m+n}R_{\alpha c\alpha c}(x_{0}).\nonumber
\end{align}

Through a series of complicated calculations, we obtain
\begin{align}
&\int_{S^{*}(M\times N)}{\rm
tr}\big[\sigma_{-2\overline{m}}\big((\triangle^{_{\varepsilon}M\times_{f} N, g^{\varepsilon, f}})\circ(\triangle^{M\times N, g})^{-\overline{m}}\big)\big](x_{0}, \xi)\nonumber\\
&=(\frac{m}{12\varepsilon}+\frac{n}{12f^{2}})area(S_{2\overline{m}})\sum_{\gamma,j=1}^{2\overline{m}}R_{\gamma j\gamma j}(x_{0})-\frac{1}{6\varepsilon}area(S_{2\overline{m}})\sum_{\gamma=1}^{2\overline{m}}\sum_{a=1}^{m}R_{\gamma a\gamma a}(x_{0})\nonumber\\
&-\frac{1}{6f^{2}}area(S_{2\overline{m}})\sum_{\gamma=1}^{2\overline{m}}\sum_{c=m+1}^{m+n}R_{\gamma c\gamma c}(x_{0}).\nonumber
\end{align}

It is easy to check that 
\begin{align}
\sum_{\gamma=m+1}^{m+n}\sum_{a=1}^{m}R_{\gamma a\gamma a}(x_{0})=\sum_{\gamma=m+1}^{m+n}\sum_{a=1}^{m}\langle R(e_{\gamma}, e_{a})e_{\gamma}, e_{a}\rangle(x_{0})=0,\nonumber
\end{align}
for this reason
\begin{align}
\sum_{\gamma=1}^{2\overline{m}}\sum_{a=1}^{m}R_{\gamma a\gamma a}(x_{0})=\sum_{\gamma,a=1}^{m}R_{\gamma a\gamma a}(x_{0}).\nonumber
\end{align}
Similarly,
\begin{align}
\sum_{\gamma=1}^{2\overline{m}}\sum_{c=m+1}^{m+n}R_{\gamma c\gamma c}(x_{0})=\sum_{\gamma,c=m+1}^{m+n}R_{\gamma c\gamma c}(x_{0}).\nonumber
\end{align}

In combination with the calculation,
\begin{align}
&\int_{S^{*}(M\times N)}{\rm
tr}\big[\sigma_{-2\overline{m}}\big((\triangle^{_{\varepsilon}M\times_{f} N, g^{\varepsilon, f}})\circ(\triangle^{M\times N, g})^{-\overline{m}}\big)\big](x_{0}, \xi)\nonumber\\
&=\big(\frac{m-2}{12\varepsilon}+\frac{n}{12f^{2}}\big)area(S_{2\overline{m}})\sum_{\gamma,a=1}^{m}R_{\gamma a\gamma a}(x_{0})+\big(\frac{m}{12\varepsilon}+\frac{n-2}{12f^{2}}\big)area(S_{2\overline{m}})\sum_{\gamma,c=m+1}^{m+n}R_{\gamma c\gamma c}(x_{0}).\nonumber
\end{align}

Motivated by Theorem 1.1, we define
\begin{defn}
Let $(A, H, D_{1})$ and $(A, H, D_{2})$ be two $2m$-dimensional spetral triple, we define the double spectral Einstein-Hilbert action as follows:
\begin{align}
EH(D_{1}, D_{2})=Wres(D_{1}^{2}D_{2}^{-2m}).\nonumber
\end{align}
\end{defn}
When $D_{1}=D_{2},$ we get the spectral Einstein-Hilbert action.

\begin{defn}
Let $g_{2}$ be Riemannian metric and $g_{1}$ is any metric (including the Lorentz metric), we define the bimetric spectral Einstein-Hilbert action as follows:
\begin{align}
EH(\triangle^{g_{1}}, \triangle^{g_{2}})=Wres\big(\triangle^{g_{1}}\circ(\triangle^{g_{2}})^{-m}\big).\nonumber
\end{align}
\end{defn}

\section{ Example }
Set $M=S^{1},$ $N=S^{3},$ where $S^{1}$ is a circle, $S^{3}$ is a there-dimensional sphere.
Then the warped product $_{\varepsilon}S^{1}\times_{f} S^{3}$ with $\varepsilon\neq0$ and a smooth function $f: S^{1}\rightarrow (0, \infty)$ for which $f>0$ is a product manifold of form $S^{1}\times S^{3}$ with the metric tensor $g^{\varepsilon, f}=\varepsilon g^{S^{1}}\oplus f^{2}g^{S^{3}}.$ 
Here $g^{S^{1}}$ is the standard Riemannian metric on $S^{1}$ inherited from $R^{2},$ $g^{S^{3}}$ is the standard Riemannian metric on $S^{3}$ inherited from $R^{4}.$
We call $(S^{1}\times S^{3}, g^{\varepsilon, f})$ the Robertson-Walker space.

Theorem 1.1 now leads to
\begin{align}
&Wres [(\triangle^{_{\varepsilon}S^{1}\times_{f} S^{3}, g^{\varepsilon, f}})\circ(\triangle^{S^{1}\times S^{3}, g})^{-2}]\nonumber\\
&=\int_{S^{1}\times S^{3}}\frac{2\pi^{2}}{\Gamma(2)}\Big[\big(\frac{1-2}{12\varepsilon}+\frac{3}{12f^{2}}\big)S_{S^{1}}+\big(\frac{1}{12\varepsilon}+\frac{3-2}{12f^{2}}\big)S_{S^{3}}\Big]d{\rm Vol_{S^{1}} }d{\rm Vol_{ S^{3}} }.\nonumber
\end{align}
Obviously, $S_{S^{1}}=0,$ $S_{S^{3}}=6.$
It follows that
\begin{align}
Wres [(\triangle^{_{\varepsilon}S^{1}\times_{f} S^{3}, g^{\varepsilon, f}})\circ(\triangle^{S^{1}\times S^{3}, g})^{-2}]&=12\pi^{2}\int_{S^{1}\times S^{3}}\big(\frac{1}{12\varepsilon}+\frac{1}{12f^{2}}\big)d{\rm Vol_{ S^{1}} }d{\rm Vol_{ S^{3}} }\nonumber\\
&=4\pi^{4}\int_{S^{1}}\big(\frac{1}{\varepsilon}+\frac{1}{f^{2}}\big)d{\rm Vol_{S^{1}}}.\nonumber
\end{align}

\section{ Declarations }

Ethics approval and consent to participate No applicable.\\

Consent for publication No applicable.\\

Availability of data and material The authors confirm that the data supporting the findings of this study are available within the article.\\

Competing interests The authors declare no conflict of interest.\\

Funding This research was funded by National Natural Science Foundation of China: No.11771070.\\

Authors' contributions All authors contributed to the study conception and design. Material preparation, data collection and analysis were performed by Siyao Liu and Yong Wang. The first draft of the manuscript was written by Siyao Liu and all authors commented on previous versions of the manuscript. All authors read and approved the final manuscript.\\

\vskip 1 true cm


\bigskip
\bigskip

\noindent {\footnotesize {\it S. Liu} \\
{School of Mathematics and Statistics, Changchun University of Science and Technology, Changchun 130022, China}\\
{Email: liusy719@nenu.edu.cn}

\noindent {\footnotesize {\it Y. Wang} \\
{School of Mathematics and Statistics, Northeast Normal University, Changchun 130024, China}\\
{Email: wangy581@nenu.edu.cn}

\end{document}